\documentclass[12pt]{amsart}
\usepackage{amsmath,amssymb,amsthm,mathtools}
\usepackage{mathrsfs,bm}
\usepackage{thmtools,thm-restate}
\usepackage{tikz-cd}
\usepackage{tikz}
\usepackage{enumitem}
\usepackage{graphicx}
\usetikzlibrary{matrix,arrows.meta,bending}
\usepackage{color}
\usepackage[numbers]{natbib}
\usepackage{booktabs}
\usepackage{longtable}
\usepackage{pgfplots}
\pgfplotsset{compat=1.14} 
\usetikzlibrary{arrows.meta,decorations.pathreplacing,positioning,calc} 

\usepackage[a4paper, left=3cm, right=3cm, top=3cm, bottom=3cm, footskip=15mm]{geometry}

\newtheorem{theorem}{Theorem}[section]

\newtheorem{lemma}[theorem]{Lemma}
\newtheorem{proposition}[theorem]{Proposition}

\theoremstyle{definition}
\newtheorem{definition}[theorem]{Definition}

\theoremstyle{remark}
\newtheorem{remark}[theorem]{Remark}

\numberwithin{equation}{section}

\makeindex

\usepackage{fancyhdr}
\usepackage{calc} 

\AtBeginDocument{%
  \setlength{\textwidth}{\dimexpr\paperwidth - 6cm\relax}%
  \setlength{\oddsidemargin}{\dimexpr 3cm - 1in\relax}%
  \setlength{\evensidemargin}{\dimexpr 3cm - 1in\relax}%
  \setlength{\marginparwidth}{2.5cm}%
}

\begin{document}

\title{On a variant of Brocard's problem via the diagonalization method
}

\author{Theophilus Agama}
\address{Department of Mathematics, African institute for mathematical sciences, Ghana,
Accra}
\email{Theophilus@ims.edu.gh/emperordagama@yahoo.com}



\date{\today}



\begin{abstract}
In this paper, we introduce and develop the method of diagonalization of functions $f:\mathbb{N}\longrightarrow \mathbb{R}$. We apply this method to show that equations of the form $\Gamma_r(n)+k=m^2$ have a finite number of solutions $n\in \mathbb{N}$ with $n>r$ for a fixed $k,r\in \mathbb{N}$, where $\Gamma_r(n)=n(n-1)\cdots (n-r)$ is the $r^{th}$ truncated Gamma function.
\end{abstract}

\maketitle

\section{Introduction and problem statement}

The classical Brocard problem - asking for integer solutions to
$$
n!+1=m^2
$$
- is a succinct and long-standing Diophantine question dating back to the work of {\it Brocard} and later noted by {\it Ramanujan}. The only small solutions known are $n=4,5,7$, and no further solutions have been found in computational searches to large bounds; see \cite{berndt2000brocard} for a concise account and computational remarks. Heuristics and extensive computation suggest that additional solutions, if any, are rare; indeed, conditional results and extensions have been established in the literature: Overholt demonstrated the finiteness of solutions under the {\sc abc} conjecture \cite{overholt1993diophantine}, while subsequent authors have considered shifted and polynomial generalizations such as $n! + A = y^2$ and equations of the form $n! = P(x)$ for polynomials $P$ of degree at least two \cite{dabrowski1996diophantine, luca2002diophantine}.\\

Motivated by these threads, the present paper studies a natural family of Brocard-type equations obtained by replacing the factorial by a truncated product. For a fixed integer $r\geq 0$ define the \emph{$r$-th truncated Gamma function}
$$
\Gamma_r(n):=\begin{cases}n(n-1)\cdots(n-r), & n>r,\\[4pt]
0, & n\leq r.\end{cases}
$$
We investigate the Diophantine equation
$$
\Gamma_r(n)+k=m^2,
$$
for fixed $k,r\in\mathbb{N}$, and we prove that for each fixed pair $(r,k)$ there are only finitely many integer solutions $n>r$. This can be viewed as a variational Brocard statement: the truncated product $\Gamma_r$ is a natural finite analogue of $n!$, and its algebraic and analytic structure admits an approach that differs from those previously used for the factorial.
\bigskip

\subsection*{Key idea and method}

The core contribution of this paper is the systematic deployment of the \emph{diagonalization method} for integer-valued sequences $f:\mathbb{N}\to\mathbb{R}$. Informally, for a fixed shift $k$, we consider the \emph{diagonal set}
$$
\mathcal{D}_k(f):=\{n\in\mathbb{N}~:~f(n)+k \text{ is a perfect square}\},
$$
and study the distribution and cumulative behaviour of \(f\) restricted to that set. Two analytic ingredients are central:\\

\begin{enumerate}
\item \textbf{Trace and integration-by-parts.} By introducing the $s$-level trace
$$
\mathbb{T}_f(s,k):=\sum_{\substack{n\leq s\\ n\in\mathcal{D}_k(f)}} f(n),
$$
and expressing it via a Stieltjes integral, we obtain relations between pointwise counts $|\mathcal{D}_k(f,s)|$ and aggregate quantities such as \(\sum_{n\le s} f(n)\) and certain Sobolev-type norms of \(f\). This reduces the discrete counting problem to estimates accessible by classical inequalities.
\bigskip

\item \textbf{A diagonal inequality via Cauchy--Schwarz.} Applying Cauchy--Schwarz to the resulting integrals produces a flexible inequality (Theorem \ref{inequality 1}) that bounds the $L^2$-trace of the diagonal in terms of normalized averages of \(f\) and a control term coming from $||f'||_{L^2}$. The inequality isolates a simple-to-verify sufficient condition which implies that $|\mathcal{D}_k(f)|$ is finite.
\end{enumerate}
\bigskip

When the general diagonal framework is applied to $f=\Gamma_r$ one benefits from explicit algebraic growth rates: $\Gamma_r(n)\asymp n^{r+1}$, and one can directly estimate the derivative-type term arising in the diagonal inequality. Elementary comparison and summation arguments then yield the finiteness result for the truncated Gamma equation; the argument is unconditional and does not rely on conjectures such as {\sc abc}.

\subsection*{Relation to previous work}

The literature on Brocard-type equations is rich: besides computational verifications and classical references (see \cite{berndt2000brocard}), conditional finiteness results based on deep conjectures appear in \cite{overholt1993diophantine}, and shifted variants and polynomial generalizations are treated in \cite{dabrowski1996diophantine, luca2002diophantine}. The diagonalization method introduced here differs in spirit: rather than embedding factorial expressions into algebraic or modular frameworks or invoking conjectural height bounds, we use a variational-analytic viewpoint that converts the square-shift condition into inequalities for aggregated traces. This viewpoint is amenable to other polynomially growing arithmetic functions and suggests a toolkit to further explore Brocard-type equations where the growth and smoothness of $f$ can be controlled.

\subsection*{Structure of the paper}

The paper is organized as follows. In Section 2, we formalize the diagonalization language: definitions (diagonals, traces, diagonal squares), the Stieltjes representation of the trace, and the main diagonal inequality (Theorem \ref{inequality 1}). Section 3 develops auxiliary estimates and the lemmas controlling the relevant $L^2$ and summatory quantities; these include Lemmas \ref{control 1} and \ref{control 2}, which quantify the behaviour of truncated products and their derivatives. Section 4 applies the diagonal method to the truncated Gamma function and proves the main theorem (Theorem \ref{Brocard}), showing finiteness of solutions to $\Gamma_r(n)+k=m^2$ for fixed $r,k$. Finally, Section 5 contains remarks, possible extensions, and examples: we indicate how the method can be adapted to other polynomially growing sequences and discuss limitations and natural directions for future work.
\medskip

The remainder of the paper develops the \emph{diagonalization} method in detail and show that the equation $\Gamma_r(n)+k=n^2$ has a finite number of solutions for $n\in \mathbb{N}$.. Our main unconditional result is stated and proved in Theorem \ref{Brocard}.

\subsection{Notations}

In this paper, we will write $f(n)\ll g(n)$ to mean there exists an absolute constant $c>0$ such that for all sufficiently large $n$, then $f(n)\leq c|g(n)|$. Conversely, we will write $f(n)\gg g(n)$ if the reverse inequality holds for all sufficiently large values of $n$. If both inequalities hold then we write in simple terms $f(n)\asymp g(n)$.

\section{The notion of diagonalization}

In this section we introduce and study the notion of \emph{diagonalization} of a function. We study this notion together with associated statistics and explore some applications.

\begin{definition}
Let $f:\mathbb{N}\longrightarrow \mathbb{R}$. Then we say $f$ is $k$-step \emph{diagonalizable} at the \emph{spot} $n\in \mathbb{N}$ if there exists some $m\in \mathbb{N}$ such that 
\begin{align}
f(n)+k=m^2.\nonumber
\end{align}
We call the set of all spots $n\in \mathbb{N}$ such that $f$ is $k$-step diagonalizable the $k^{th}$-step diagonal of $f$ and denote by $\mathcal{D}_k(f)$. We call the set of all truncated spots $\mathcal{D}_k(f)\cap \mathbb{N}_s:=\mathcal{D}_k(f,s)$ the $s^{th}$ scale diagonal. We call the set of all squares 
\begin{align}
\mathbb{B}_k(f):=\left \{m^2\in \mathbb{N}~|~f(n)+k=m^2\right \}\nonumber
\end{align}
the $k^{th}$-step diagonal squares. We write the length of this diagonal as 
\begin{align}
|\mathcal{D}_k(f,s)|:=\# \{n\leq s~|~f(n)+k=m^2\}.\nonumber
\end{align}
\end{definition}
\bigskip

It is easy to see that $|\mathcal{D}_k(f,s)|<s$.

\subsection{The $s$-level trace of the diagonal}

In this section, we introduce the notion of the \emph{trace} of the diagonal. We launch and examine the following languages.

\begin{definition}
By the $s^{th}$ level trace of the diagonal $\mathcal{D}_k(f)$, denoted $\mathbb{T}_{f}(s,k)$, we mean the partial sum 
\begin{align}
\mathbb{T}_{f}(s,k):=\sum \limits_{\substack{n\leq s\\n\in \mathcal{D}_k(f)}}f(n).\nonumber
\end{align}
\end{definition}
\bigskip

Let us suppose that $f$ is a function with continuous derivative on $[1,s]$ for $s\geq 1$ with $s\in \mathbb{R}$, then applying the Stieltjes integration by parts, we can write the $s^{th}$ level trace of the diagonal in the form 

\begin{align}
\mathbb{T}_{f}(s,k):&=\sum \limits_{\substack{n\leq s\\n \in \mathcal{D}_k(f)}}f(n)\nonumber \\&=\int \limits_{1^{-}}^{s}f(t)d|\mathcal{D}_k(f,t)|\nonumber \\&=f(s)|\mathcal{D}_k(f,s)|-\int \limits_{1}^{s}f'(t)|\mathcal{D}_k(f,t)|dt.\nonumber
\end{align}

\begin{theorem}[Diagonal inequality]\label{inequality 1}
Let $f$ be a function with continuous derivative on $[1,s]$ for $s\geq 1$ with $s\in \mathbb{R}$. If 
\begin{align}
\mathcal{D}_k(f,s)-\frac{1}{f(s)}\bigg(\sqrt{\int \limits_{1}^{s}|f'(t)|^2dt}\bigg)\times \bigg(\int \limits_{1}^{s}|\mathcal{D}_k(f,t)|^2dt\bigg)^{\frac{1}{2}}\geq 0\nonumber
\end{align}
for all $s\geq 1$, then the inequality
\begin{align}
(\int \limits_{1}^{s}|\mathcal{D}_k(f,t)|^2dt)^{\frac{1}{2}} \ll \bigg(\frac{1}{f(s)}\sum \limits_{n\leq s}f(n)\bigg)\bigg(1-\frac{1}{f(s)}\sqrt{\int \limits_{1}^{s}|f'(t)|^2dt}\bigg)^{-1}\nonumber
\end{align}
holds.
\end{theorem}

\begin{proof}
Following the ensuing discussion, we obtain the upper bound
\begin{align}
|\mathcal{D}_k(f,s)|&\leq \frac{1}{f(s)}\sum \limits_{n\leq s}f(n)+\frac{1}{f(s)}\int \limits_{1}^{s}f'(t)|\mathcal{D}_k(f,t)|dt\nonumber
\end{align}
so that by using the Cauchy-Schwartz inequality, we further obtain the upper bound
\begin{align}
|\mathcal{D}_k(f,s)|&\leq \frac{1}{f(s)}\sum \limits_{n\leq s}f(n)+\frac{1}{f(s)}\bigg(\int \limits_{1}^{s}|f'(t)|^2dt\bigg)^{\frac{1}{2}}\times \bigg(\int \limits_{1}^{s}|\mathcal{D}_k(f,t)|^2dt\bigg)^{\frac{1}{2}}.\nonumber
\end{align}
Rearranging terms, using the condition
\begin{align}
\mathcal{D}_k(f,s)-\frac{1}{f(s)}\bigg(\sqrt{\int \limits_{1}^{s}|f'(t)|^2dt}\bigg)\times \bigg(\int \limits_{1}^{s}|\mathcal{D}_k(f,t)|^2dt\bigg)^{\frac{1}{2}}\geq 0 \nonumber
\end{align} 
and noting that 
\begin{align}
(\int \limits_{1}^{s}|\mathcal{D}_k(f,t)|^2dt)^{\frac{1}{2}}\geq |\mathcal{D}_k(f,s)|\nonumber
\end{align}
for all $s\geq 1$, then the claimed inequality holds.
\end{proof}
\bigskip

\section{The diagonal method}

Brocard's problem asks if there are a finite number of solutions to equation $n!+1=m^2$. Proposition \ref{inequality 1} provides a helpful inequality to explore. The current framework can be used to investigate a much broader version of the problem. We may improve the outcome by using the \textbf{Diagonal} inequality.

\begin{proposition}[The Diagonal method]\label{general Brocard-trace}
Let
\begin{align}
|\mathcal{D}_k(f,s)|-\frac{1}{f(s)}\bigg(\sqrt{\int \limits_{1}^{s}|f'(t)|^2dt}\bigg)\times \bigg(\int \limits_{1}^{s}|\mathcal{D}_k(f,t)|^2dt\bigg)^{\frac{1}{2}}\geq 0\nonumber
\end{align}
for all $s\geq 1$. If 
\begin{align}
\lim \limits_{s\longrightarrow \infty}\bigg(\frac{1}{f(s)}\sum \limits_{n\leq s}f(n)\bigg)\bigg(1-\frac{1}{f(s)}\sqrt{\int \limits_{1}^{s}|f'(t)|^2dt}\bigg)^{-1}<\infty \nonumber
\end{align}
then equation $f(n)+k=m^2$ has only a finite number of solutions in $\mathbb{N}$ for a fixed $k\in \mathbb{N}$.
\end{proposition}

\begin{proof}
Using Proposition \ref{inequality 1}, it follows under the requirements that $|\mathcal{D}_k(f)|<\infty$, and the claim follows immediately.
\end{proof}

\begin{remark}
The upper bound derived from Proposition \ref{inequality 1} provides a useful tool for studying the size of the quantity 
$$
\# \{n\leq s~|~f(n)+k=m^2\}
$$ 
and in particular Brocard's problem, which asks if the set of integers whose factorials are unit left-translate of a square is either an infinite set or a finite set. It is worth noting that the upper bounds we have derived do not depend on the size of the shift, but on the underlying function. This uniformity suggests the actual size of the quantity 

\begin{align}
\# \{n\leq s~|~f(n)+k=m^2\}\nonumber
\end{align}
will be mainly influenced by the function under consideration. In various circumstances, the ease with which to verify the underlying conditions will inform the category of bounds to exploit. Now, we apply the \emph{Diagonal method} to study a slight variant of Brocard's problem.
\end{remark}

\begin{lemma}\label{control 1}
We have
\begin{align}
\bigg(\int \limits_{1}^{s}|\mathcal{D}_k(f,t)|^2dt\bigg)^{\frac{1}{2}} \ll |\mathcal{D}_k(f,s)|^{\frac{3}{2}}.\nonumber
\end{align}
\end{lemma}

\begin{remark}
The upper bound of the lemma \ref{control 1} can easily be obtained by exploiting the techniques of integrating a function in elementary calculus.
\end{remark}
\bigskip

\begin{definition}[The $r^{th}$ truncated Gamma function]
Let $r\in \mathbb{N}$ be fixed. By the $r^{th}$ \emph{truncated} Gamma function $\Gamma_r$, we mean the function 
\begin{align}
\Gamma_r(n):=
\begin{cases}n(n-1)\cdots (n-r)~\text{if}~n>r\\0~\text{otherwise.}
\end{cases}\nonumber
\end{align} 
\end{definition}
\bigskip

It follows from the definition of the $r^{th}$ truncated Gamma function that for $s>r$ then
$$
\Gamma_r(s)=s(s-1)\cdots (s-r).
$$
It follows in a similar way that 
$$
\Gamma_r(s-1)=(s-1)(s-2)\cdots (s-1-(r-1))
$$ 
since $s-1>r-1$ and 
$$
\Gamma_r(s-2)=(s-2)\cdots (s-2-(r-2))
$$ 
since $s-2>r-2$ and so on. The $r^{th}$ truncated Gamma function for smaller arguments can be obtained in a similar manner.

\begin{lemma}\label{control 2}
For all $s>r$, we have 
\begin{align}
\frac{1}{\Gamma_r(s)}\sqrt{\int \limits_{1}^{s}|\Gamma_r'(t)|^2dt} \asymp \frac{1}{s^{\frac{1}{2}}}.\nonumber
\end{align}
\end{lemma}

\begin{proof}
It follows naturally from the definition of the $r^{th}$ \textit{truncated} Gamma function that $\Gamma_r(s)\asymp s^{r+1}$ so that 
\begin{align}
\int \limits_{1}^{s}|\Gamma_r'(t)|^2dt \asymp s^{2r+1}\nonumber
\end{align}
and the claimed upper bound is an easy consequence.
\end{proof}

\section{An application}

\begin{theorem}[Variational Brocard] \label{Brocard}
The equation $\Gamma_r(s)+k=m^2$ has a finite number of solutions $s\in \mathbb{N}$ with $s>r$ for a fixed $k,r\in \mathbb{N}$, where $\Gamma_r$ is the $r^{th}$ truncated Euler Gamma function.
\end{theorem}

\begin{proof}
We first apply the lemma \ref{control 2} and note that 
\begin{align}
\bigg(\int \limits_{1}^{s}|\mathcal{D}_k(\Gamma_r,t)|^2dt\bigg)^{\frac{1}{2}} \ll |\mathcal{D}_k(\Gamma_r,s)|^{\frac{3}{2}} \ll  |\mathcal{D}_k(\Gamma_r,s)|\sqrt{s} \nonumber
\end{align}
since $|\mathcal{D}_k(\Gamma_r,s)|<s$. It suffices to check 
\begin{align}
\lim \limits_{s\longrightarrow \infty}\frac{1}{\Gamma_r(s)}\sum \limits_{n\leq s}\Gamma_r(n)<\infty \nonumber
\end{align}
and
\begin{align}
\lim \limits_{s\longrightarrow \infty}\frac{1}{\Gamma_r(s)}\sqrt{\int \limits_{1}^{s}|\Gamma_r'(t)|^2dt}<\infty \nonumber
\end{align}
and using the inequality 
$$
(\int \limits_{1}^{s}|\mathcal{D}_k(\Gamma_r,t)|^2dt)^{\frac{1}{2}} \geq |\mathcal{D}_k(\Gamma_r,s)|.
$$
\end{proof}

\section{Further remarks}

The \emph{diagonal} method developed can be used to investigate the original problem of Brocard concerning the finiteness of positive integer solutions to the equation $n!+1=m^2$. The Brocard problem is equivalent to asking whether the equation $\Gamma(n+1)+1=m^2$ has a finite integer solutions $n,m$. In its general form, one may consider the equation with general shifts $\Gamma(n+1)+k=m^2$ for a fixed $k\geq 1$ with $k\in \mathbb{N}$. It will suffice to verify the requirements in the diagonal method in the case $f:=\Gamma$ to address the Brocard problem.  

\bibliographystyle{amsplain}

\end{document}